\documentclass{amsart}
\usepackage{amsthm}
\usepackage{amssymb}
\usepackage{amsmath}
\usepackage{caption}
\usepackage{graphicx}
\usepackage{mathrsfs}
\usepackage{lscape}
\usepackage{float}
\usepackage{color}
\usepackage{amscd}
\usepackage{booktabs}
\usepackage{hyperref}

\newcommand{\Z}{\mathbb{Z}}
\newcommand{\F}{\mathbb{F}}
\newcommand{\N}{\mathbb{N}}
\newcommand{\R}{\mathbb{R}}
\newcommand{\E}{\mathbb{E}}
\newcommand{\Q}{\mathbb{Q}}

\newcommand{\T}{\mathcal{T}}
\newcommand{\Y}{\mathcal{Y}}
\newcommand{\expect}{\mathbb{E}}
\DeclareMathOperator{\Aut}{Aut}

\newcommand{\Erdos}{Erd\H{o}s }
\newcommand{\Renyi}{R\'{e}nyi }
\usepackage{relsize}
\newtheorem{theorem}{Theorem}

\theoremstyle{definition}
\newtheorem{conjecture}{Conjecture}
\newtheorem{definition}{Definition}

\newtheorem{example}{Example}

\title[Cohen--Lenstra heuristics for torsion in homology]{Cohen--Lenstra heuristics for torsion in homology of random complexes}
\author{Matthew Kahle}
\address{The Ohio State University}
\email{mkahle@math.osu.edu}
\thanks{partially supported by NSF-DMS \#1352386}
\author{Frank H.~Lutz}
\address{Technische Universit\"at Berlin}
\email{lutz@math.tu-berlin.de}
\thanks{partially supported by the DFG Coll. Research Center TRR 109, ``Discretization in Geometry and Dynamics''}
\author{Andrew Newman}
\address{The Ohio State University}
\email{newman.534@osu.edu}
\thanks{partially supported by NSF-DMS \#1547357}
\author{Kyle Parsons}
\address{The Ohio State University}
\email{parsons.299@osu.edu}
\thanks{partially supported by NSF-DMS \#1547357}
\date{\today}

\begin{document}

\maketitle

\begin{abstract}
We study torsion in homology of the random $d$-complex $Y \sim Y_d(n,p)$ experimentally. Our experiments suggest that there is almost always a moment in the process where there is an enormous burst of torsion in homology $H_{d-1}(Y)$. This moment seems to coincide with the phase transition studied in \cite{AL,LP,LP3} , where cycles in $H_d(Y)$ first appear with high probability.

Our main study is the limiting distribution on the $q$-part of the torsion subgroup of $H_{d-1}(Y)$ for small primes $q$. We find strong evidence for a limiting Cohen--Lenstra distribution, where the probability that the $q$-part is isomorphic to a given $q$-group $H$ is inversely proportional to the order of the automorphism group $|\mbox{Aut}(H)|$.

We also study the torsion in homology of the uniform random $\Q$-acyclic $2$-complex. This model is analogous to a uniform spanning tree on a complete graph, but more complicated topologically since Kalai showed that the expected order of the torsion group is exponentially large in $n^2$  \cite{Kalai}. We give experimental evidence that in this model also, the torsion is Cohen--Lenstra distributed in the limit.
\end{abstract}

\section{Introduction}

The random $d$-complex $Y_d(n, p)$, introduced by Linial and Meshulam in \cite{LM}, is the probability space 
on $d$-dimensional simplicial complexes with vertex set $[n]$ and complete $(d-1)$-skeleton, where each $d$-dimensional face is included independently with probability $p$. A closely related model is $Y_d(n,m)$, which also has vertex set $[n]$ 
and complete $(d-1)$-skeleton, and has exactly $m$ $d$-dimensional faces, chosen uniformly at random 
from all $\left(\!\!\genfrac{}{}{0pt}{}{\binom{n}{d+1}}{m}\!\!\right)$ possibilities.

The random complex $Y_d(n,m)$ is a generalization of the \Erdos--\Renyi random graph $G(n,m)$. By now, many topological properties 
of the random $d$-complex have been studied \cite{AL, AL2, ALLM, HKP, LM, LP, MW, LP2}. Many of these properties 
are \emph{monotone} properties which always have thresholds. Torsion in homology is a  non-monotone property, though, 
so even formulating well-posed questions requires a bit more care.  

We regard the random $d$-complex as a discrete-time stochastic process, as follows. At time $m=0$, we start with the $(d-1)$-skeleton 
of the simplex on $n$ vertices, and for each time $m=1, 2, \dots, {n \choose d+1}$, we add one face $\sigma$, chosen uniformly at random 
from all remaining $d$-dimensional faces. Following the notation of \cite{LP2}, we let $\Y_d(n)$ denote the discrete-time 
stochastic process $\Y_d(n) = \{ Y_d(n,m) :  0 \le m \le {n \choose d+1} \}$.

We always assume that homology is with integer coefficients, unless otherwise indicated. Our main interest is in the torsion in $H_{d-1}(\Y_d(n))$. Experimentally, there is a short burst of huge torsion---an example run is shown in Table \ref{tbl:75vertices}.
\begin{table}[H]
\small\centering
\defaultaddspace=0.2em
\caption{The torsion burst for a single instance of $\mathcal{Y}_2(75)$.}\label{tbl:75vertices}
\begin{tabular}{@{}r@{\hspace{12mm}}l@{\hspace{12mm}}l@{}}
\\\toprule
 \addlinespace
 \addlinespace
    $m$  & $H_2$  & $H_1$    \\ \midrule
 \addlinespace
 \addlinespace
 \addlinespace
2469 & $\Z^4$ & $\Z^{236}$ \\
2470 & $\Z^4$ & $\Z^{235} \times \Z/2\Z$ \\
2471 & $\Z^4$ & $\Z^{234} \times \Z/2\Z$ \\
2472 & $\Z^4$ & $\Z^{233} \times \Z/2\Z$ \\
2473 & $\Z^4$ & $\Z^{232} \times \Z/2\Z$ \\
2474 & $\Z^4$ & $\Z^{231} \times \Z/2\Z \times \Z/2\Z$ \\
2475 & $\Z^4$ & $\Z^{230} \times \Z/2\Z \times \Z/79040679454167077902597570\Z$ \\
2476 & $\Z^5$ & $\Z^{230} \times \Z/2\Z$\\
2477 & $\Z^5$ & $\Z^{229} \times \Z/2\Z$\\
2478 & $\Z^6$ & $\Z^{229}$\\
 \addlinespace
 \addlinespace
\bottomrule
\end{tabular}
\end{table}


The torsion burst seems closely related to the appearance of the first nontrivial top-dimensional homology class, that is, the first homology class
not represented by an embedded copy of $\partial \Delta^{d + 1}$. 
In \cite{AL} and \cite{LP} a constant $c_d$ is described such that $p = c_d/n$ is the sharp threshold for $H_d(Y_d(n, p), \R) \neq 0$.
The connection between the emergence of nontrivial top homology and torsion in the $(d-1)$st homology was mentioned in \cite{LP2} 
where they make the following conjecture that there is no torsion away from the phase transition. Additionally, \cite{LP2} 
mentions conducting experiments indicating torsion near the phase transition at $c_d/n$. The torsion burst was also observed 
for $d = 2$ in experiments examining the fundamental group of random 2-complexes in \cite{BL}.

\begin{conjecture}[\L uczak and Peled  \cite{LP2}]
For every $d \geq 2$ and $p = p(n)$ such that $|np - c_d|$ is bounded away from 0, $H_{d-1}(Y_d(n, p); \Z)$ is torsion-free with high probability. 
\end{conjecture} 

It may seem surprising that there is a simplicial complex on 75 vertices which has the huge torsion we see in Table \ref{tbl:75vertices}, and we don't yet know why this torsion appears, however the existence of small simplicial complexes with large torsion in homology was known previously.  Indeed, there is a canonical class of examples which realizes this phenomenon. 
These complexes were first described by Kalai in \cite{Kalai} and are called $\Q$-acyclic complexes.  A $d$-dimensional $\Q$-acyclic complex 
is defined to be a $d$-dimensional simplicial complex $X$ with complete $(d-1)$-skeleton, $\binom{n-1}{d}$ $d$-dimensional faces, 
$\beta_d(X, \mathbb{Q}) = 0$, and $\beta_{d - 1}(X, \Q) = 0$. Thus, $\Q$-acyclic complexes are higher dimensional generalizations of trees. 
However, unlike trees, $d$-dimensional $\Q$-acyclic complexes may have finite but nontrivial $(d - 1)$st homology group. 

In \cite{Kalai}, Kalai gives a striking generalization of Cayley's enumeration of spanning trees, showing 
$$\sum_{X \in \T^d(n)} |H_{d - 1}(X)|^2 = n^{\binom{n - 2}{d}},$$
where $\T^d(n)$ is defined to be the collection of all $d$-dimensional $\Q$-acyclic complexes on $n$ vertices. 
As a corollary to the main result, \cite{Kalai} shows that 
$$\expect \left[ |H_{d - 1}(X)| \right]\geq \exp \left( \Theta \left(n^d \right) \right)$$
for $X$ a uniform random $\Q$-acyclic complex with $d \geq 2$, and also points out that $|H_{d - 1}(X)| \leq \sqrt{d + 1}^{\binom{n-2}{d}}$ 
for any simplicial complex $X$ having dimension $d$ on $n$ vertices. So the existence of small complexes with large torsion 
in homology is known, however, the reason for them to appear in $\Y_2(n)$ remains a mystery.
%

\section{Background on Cohen--Lenstra heuristics}

Cohen--Lenstra heuristics refer to a family of probability distributions on finite abelian groups, first appearing in \cite{CL} in the context of class groups of number fields, in which the probability of each  
group $G$ in the support of the distribution is inversely proportional to $|G|^{\alpha}|\Aut(G)|^{\beta}$, 
for  $\alpha \geq 0$ and $\beta \geq 0$. Since Cohen--Lenstra heuristics were first introduced they have appeared in number theory as well as in various models of random integral matrix cokernels (see for example \cite{CKLPW, CL, EVW, FW, Koplewitz, Maples, MatchettWood2, MatchettWood3, MatchettWood}). The variety of settings in which Cohen--Lenstra heuristics appear suggest they are a natural family of distributions on finite abelian groups.

A Cohen--Lenstra heuristic of particular interest is the case where $\alpha = 0$ and $\beta = 1$. However, this does not give a distribution on the set of \emph{all} finite abelian groups. That is, it is well known that there is no distribution on the set of all finite abelian groups so that for any such group $G$,
$$\Pr(G) \propto \frac{1}{|\Aut(G)|}.$$
A simple proof of this fact appears in, for example, Chapter 5 of \cite{Lengler}. However, if one restricts to the family of $p$-groups $\mathcal{G}_p$ for any fixed prime $p$ then there is such a distribution so that the probability assigned to any $p$-group is inversely proportional to the number of automorphisms of that group. Indeed, Cohen and Lenstra prove in their original paper \cite{CL} that for a fixed prime $p$,
$$\sum_{G \in \mathcal{G}_p} \frac{1}{|\Aut(G)|} = \frac{1}{\prod_{i = 1}^{\infty}(1 - p^{-k})}.$$
This family of distributions on $p$-groups is the one we consider most often here. However, in Section ~\ref{sec:othergroups} we also discuss  families of distributions with $\beta = 1$, but $\alpha$ equal to a positive integer. \\

In the next section, we focus on the $2$-dimensional process $\Y_2(n)$ and experimentally measure 
where the torsion burst occurs, how large the torsion group is at its peak, and the Cohen--Lenstra heuristics which model the random groups within the torsion burst.  
In Section~\ref{sec:higher_dim}, we discuss torsion in higher dimensions of the Linial--Meshulam process.
In Section~\ref{sec:Q-acyclic}, we experimentally measure both the size and the distribution of the first homology groups 
of uniformly random $2$-dimensional $\Q$-acyclic complexes.
In Section \ref{sec:hit}, we make more refined conjectures about the torsion burst.\\



\section{Random 2-complex experiments} \label{sec:2d}
\subsection{Preliminaries}

For $n \in \N$, we run the stochastic process $\Y_2(n)$ and compute integer homology at each step, then output the largest torsion subgroup 
that appears in $H_1$. If there are two or more nonisomorphic groups which attain the same maximum size in the process, we take the first. 
We denote this randomly-generated group $LT(n)$.

We do not actually compute $\binom{n}{3}$ rounds of integer homology computations. Instead we make a number of time-saving reductions. 
The first reduction is to consider just a portion of the entire $\binom{n}{3}$-step process.  From the comments of \L uczak and Peled 
and their conjecture \cite{LP2}, as well as our very early experiments, the torsion burst appears to occur around the top-dimensional homology threshold. 
Since \cite{AL} and \cite{LP} show that $c_2/n$ is the sharp phase transition for top-dimensional homology to emerge, we set our window 
to search for $LT(n)$ around $m^* := \lfloor \frac{c_2}{n} \binom{n}{3}\rfloor$. The exact window varied between different rounds of experiments, 
but based on our most recent trials searching the window from $m_1 = m^* - 100$ to $m_2 = m^* + 100$ appears to almost always be sufficient 
to find the torsion burst. 

Furthermore we are aided by the fact that we are looking for the largest torsion group rather than, say, the first torsion group. 
If $\beta_i$ denotes the $i$th Betti number over $\Q$, we observe that the size of the torsion group can increase only when a face is added 
which decreases $\beta_1$, and can decrease only when $\beta_2$ increases (this is because $H_1(Y_2(n, m + 1))$ is a quotient 
of $H_1(Y_2(n, m))$ by the Mayer--Vietoris sequence). Thus, if we compute $\beta_2$ with real coefficients, at each stage in the window 
that we search we only need to compute integer homology of $Y_2(n, m)$ if $\beta_2(Y_2(n, m - 1)) = \beta_2(Y_2(n, m)) < \beta_2(Y_2(n, m + 1))$. 
This seems to be a large reduction in the number of integer homology computations. Before the phase transition, it is rare 
that $\beta_2$ will increase and after the phase transition it is rare that it will not. This  is made precise by considering 
the homological shadow of $Y_2(n, p)$ described in \cite{LP} and \cite{LNPR}. We discuss the homological shadow in the final section.

Computing $\beta_1$ or $\beta_2$ with real coefficients or rational coefficients is faster than computing full integer homology, 
however, the process of computing the rank of the relevant boundary matrix typically results in the entries growing arbitrarily large, so the 
process is still slow. To save time, we instead pick a large prime $q_0$ and compute $\beta_2$ 
with $\Z/q_0\Z$ coefficients instead of computing $\beta_2$ with rational coefficients, and then we compute integer homology 
for all values of $m$ for which $\beta_2(Y_2(n, m - 1); \Z/q_0\Z) = \beta_2(Y_2(n, m); \Z/q_0\Z) < \beta_2(Y_2(n, m + 1); \Z/q_0\Z)$. 
We keep the largest torsion group found in the first homology as $LT(n)$. If there is a tie between two non-isomorphic largest torsion groups 
of the same size, we keep the first one as $LT(n)$. In practice, ties like this never occurred. 

\begin{example}
We run our implementation of $LT(n)$ with $n = 60$ and $q_0 = 10007$. First, the predicted place for torsion to appear is
 computed as $m^* = \lfloor \frac{c_2}{60} \binom{60}{3}\rfloor$. This value is 1570; so we set the window to be from $m_1 = 1470$ to $m_2 = 1670$. 
Now we generate a random $2$-complex on 60 vertices with 1670 $2$-dimensional faces which are randomly ordered. 
This random ordering gives us the first 1670 steps in an instance of $\Y_2(60)$. Now we do a binary search to find all values 
of $m$ between 1471 and 1669 so that $\beta_2(Y_2(60, m - 1); \Z/10007\Z) = \beta_2(Y_2(60, m); \Z/10007\Z) < \beta_2(Y_2(60, m + 1), \Z/10007\Z)$. 
In this particular instance there happened to be two such values of $m$, 1543 and 1545. Finally we compute integer homology 
of $Y_2(60, 1543)$ and $Y_2(60, 1545)$ and compare the torsion parts which are $\Z/66911823408\Z$ and $\Z/4\Z$, respectively. 
The output therefore is $\Z/66911823408\Z$. 
\end{example} 
%

We ran our code for these Linial--Meshulam experiments in GAP by implementing an algorithm of Benedetti and Lutz \cite{BL} 
to find a Morse matching (in the sense of discrete Morse theory \cite{Forman}). After finding a Morse matching, we build 
the top-dimensional boundary matrix for the resulting CW-complex using standard techniques and then use the built-in linear algebra functions 
of GAP to compute integer homology or homology with respect to finite field coefficients from the boundary matrix. 
 

As in other models of random abelian groups such as those studied in \cite{CL, FW, MatchettWood, MatchettWood2}, a reasonable conjecture 
is that for a fixed prime $q$, the Sylow $q$-subgroup of $LT(n)$ is distributed according to a Cohen--Lenstra distribution. Based on our data 
we make the following conjecture.
\begin{conjecture}\label{cj:LT}
For a fixed prime $q$, the Sylow $q$-subgroup of $LT(n)$ is asymptotically distributed according to the Cohen--Lenstra distribution 
which assigns probability $\dfrac{\prod_{k = 1}^{\infty} (1 - q^{-k})}{|\Aut(G)|}$ to any $q$-group~ $G$.
\end{conjecture}
To be more precise about the type of convergence, we conjecture that for any fixed prime $q$, the Sylow $q$-subgroup of $LT(n)$ converges 
in distribution to the Cohen--Lenstra heuristic given in the statement of Conjecture \ref{cj:LT}.

\subsection{Experiments and Results}

For our experiments we ran $LT(n)$ enough times so that we generated 10,000 nontrivial abelian groups. We do not think omitting trials 
where $LT(n)$ returns the trivial group changes our data too much, and consider such a situation an error state for our experiments. 
There are three reasons why $LT(n)$ might return the trivial group:
\begin{enumerate}
\item There is no torsion burst.\\[-2mm]
\item There is a torsion burst, but it occurs before or after the window that we search, so we miss it.\\[-2mm]
\item There is $q_0$-torsion which affects our search process since in this case the Betti number 
         $\beta(Y_2(n, m); \Z/q_0\Z)$ is not always equal to $\beta_2(Y_2(n, m); \Q)$. 
\end{enumerate}
Since we do not do anything to distinguish between these three possibilities we just omit and replace trials with $LT(n) = 1$. 
We do not believe this has a significant impact on our data. The simplest reason is that it is rare that $LT(n)$ returns the trivial group. 
In the trials on 50 vertices, only about 5\% returned the trivial group (that is 500 of the initial 10,000 trials returned the trivial group 
and were replaced with 500 new trials all of which returned nontrivial groups), and that the proportion dropped significantly as $n$ 
increased all the way to less than 0.1\% for trials on 125 vertices. Moreover, for our experiments here we used $q_0 = 10007$, 
and so if our conjecture is true, we only expect to miss the torsion burst due to the presence of $q_0$-torsion in about 0.01\% of trials. 
Finally, the conjectured distribution is a distribution on the $q$-part of $LT(n)$ for a fixed prime $q$. This distribution can not be extended in the naive way to a distribution on all finite abelian groups that assigns a well defined probability to every finite abelian group, but see \cite{Lengler} for a discussion of global Cohen--Lenstra heuristics. In this probability distribution, the trivial group has probability zero. Based on our conjecture and the evidence in favor of it, we expect that asymptotically almost surely, $\Y_2(n)$ has no torsion burst.

Our conjecture gives a prediction for how likely a particular $q$-group is to appear, for any particular prime $q$. In Tables \ref{tbl:2groups}, \ref{tbl:3groups}, and \ref{tbl:5groups},
 we summarize how closely our empirical distribution is to the conjectured Cohen--Lenstra distribution for $q =$ 2, 3, and 5, respectively. 
 Since the conjecture is that the probability of a group is inversely proportional to the size of its automorphism group, 
 Tables \ref{tbl:2groups}, \ref{tbl:3groups}, and \ref{tbl:5groups} show for a given group $G$ the observed ratio of the number of instances 
 where the Sylow $q$-group of $LT(n)$ was trivial to the number of instances where the Sylow $q$-group of $LT(n)$ was isomorphic to $G$. 
 This observed ratio from our experiments is compared with the conjectured ratio in Table \ref{tbl:2groups}. 
 For example, $\Z/2\Z \times \Z/2\Z$ has six automorphisms, so under the limiting conjectured distribution $LT(n)$ 
 should have trivial Sylow 2-subgroup six times as often as it has $\Z / 2\Z \times \Z/2\Z$ as its Sylow 2-subgroup.
\begin{table}
\small\centering
\defaultaddspace=0.2em
\caption{The empirical ratio of the probability that $LT(n)$ has trivial $2$-part to the probability of the given $2$-group, 
compared to the predicted ratio from the Cohen--Lenstra distribution (10,000 trials for each $n$).}\label{tbl:2groups}
\begin{tabular*}{\linewidth}{@{}l@{\extracolsep{\fill}}r@{\extracolsep{\fill}}r@{\extracolsep{\fill}}r@{\extracolsep{\fill}}r@{\extracolsep{\fill}}r@{\extracolsep{\fill}}r@{}}
\\\toprule
 \addlinespace
 \addlinespace
 $2$-Subgroup & $n = 50$ & $n = 60$ & $n = 75$ & $n = 100$ & $n = 125$ & Predicted   \\ \midrule
 \addlinespace
 \addlinespace
 \addlinespace
Trivial Group & 1 & 1 & 1 & 1 & 1 & 1\\
$\Z/2\Z$  & 0.9045 & 1.0364 & 0.9752 & 0.9965 & 1.0076 & 1 \\
$\Z/4\Z$ & 1.8822 & 1.9795 & 2.0063 & 2.0314 & 2.0423& 2 \\
$\Z/8\Z$ & 3.8530 & 4.3975 & 3.8548 & 4.2812 & 4.2525 & 4 \\
$\Z/2\Z \times \Z/2\Z$  & 5.3591 & 5.6712 & 6.0705 & 5.8846& 5.9958 & 6 \\
$\Z/2 \times \Z/4\Z$ & 6.4772 & 7.6464 & 7.4177 & 8.2820 & 7.9560 & 8 \\
$\Z/16\Z$  & 7.2026& 7.6666 & 7.1561 & 7.3040 & 7.3316 & 8 \\
$\Z/2 \times \Z/8$  & 13.1756 & 14.6363 & 14.9526 & 17.5120 & 13.7251 & 16 \\
$\Z/32\Z$ & 15.3466 & 16.8488 & 16.2342 & 19.9109 & 14.6262 & 16 \\
 \addlinespace
 \addlinespace
\bottomrule
\end{tabular*}
\end{table}


\begin{table}
\small\centering
\defaultaddspace=0.2em
\caption{The empirical ratio of the probability that $LT(n)$ has trivial $3$-part to the probability of the given $3$-group, 
compared to the predicted ratio from the Cohen--Lenstra distribution (10,000 trials for each $n$).}\label{tbl:3groups}
\begin{tabular*}{\linewidth}{@{}l@{\extracolsep{\fill}}r@{\extracolsep{\fill}}r@{\extracolsep{\fill}}r@{\extracolsep{\fill}}r@{\extracolsep{\fill}}r@{\extracolsep{\fill}}r@{}}
\\\toprule
 \addlinespace
 \addlinespace
$3$-Subgroup & $n = 50$ & $n = 60$ & $n = 75$ & $n = 100$ & $n = 125$ & Predicted   \\ \midrule
 \addlinespace
 \addlinespace
 \addlinespace
Trivial Group & 1 & 1 & 1 & 1 & 1 & 1\\
$\Z/3\Z$  & 1.9622 & 1.9277& 2.0334 & 1.9381 & 2.0482 & 2 \\
$\Z/9\Z$  & 6.0997 & 5.8486 & 5.9510 & 6.0010 & 5.9337 & 6 \\
$\Z/27\Z$ & 17.4937 & 17.8290 & 17.0212 & 17.8295 & 17.8012 & 18 \\
$\Z/3\Z \times \Z/3\Z$ & 44.1507 & 44.9349 & 47.8632 & 42.0075 & 47.8220 & 48 \\
$\Z/81\Z$  & 47.1440 & 54.7227 & 54.3689 & 56.5816 & 64.8620 & 54 \\
 \addlinespace
 \addlinespace
\bottomrule
\end{tabular*}
\end{table}


\begin{table}
\small\centering
\defaultaddspace=0.2em
\caption{The empirical ratio of the probability that $LT(n)$ has trivial $5$-part to the probability of the given $5$-group, 
compared to the predicted ratio from the Cohen--Lenstra distribution (10,000 trials for each $n$).}\label{tbl:5groups}
\begin{tabular*}{\linewidth}{@{}l@{\extracolsep{\fill}}r@{\extracolsep{\fill}}r@{\extracolsep{\fill}}r@{\extracolsep{\fill}}r@{\extracolsep{\fill}}r@{\extracolsep{\fill}}r@{}}
\\\toprule
 \addlinespace
 \addlinespace
$5$-Subgroup & $n = 50$ & $n = 60$ & $n = 75$ & $n = 100$ & $n = 125$ & Predicted   \\ \midrule
 \addlinespace
 \addlinespace
 \addlinespace
Trivial Group & 1 & 1 & 1 & 1 & 1 & 1\\
$\Z/5\Z$  & 4.1215 & 4.1488 & 4.1628 & 4.0488 & 3.9284 & 4 \\
$\Z/25\Z$ & 20.4866 & 20.7777 & 20.1259 & 19.5384 & 20.2080 & 20 \\
 \addlinespace
 \addlinespace
\bottomrule
\end{tabular*}
\end{table}


%
%

For these 10,000 trials, we also computed the total variation distances between the empirical distributions 
coming from our data sets and the conjectured distribution for primes $q=2, 3, 5, \dots, 23$. We do not include a table of these computations, since all the distances were less than $0.03$.


In addition to experimentally establishing Cohen--Lenstra heuristics, we were also interested in how large these torsion groups 
can be compared to the number of vertices in the simplicial complex. It is known from \cite{Kalai} that on average 
$d$-dimensional $\Q$-acyclic complexes have groups of size $\exp(\Theta(n^d))$ as their codimension-$1$ homology groups, 
and that the maximum size of the codimension-$1$ homology group of a $\Q$-acyclic complex is also bounded by $\exp(\Theta(n^d))$. 
So indeed the torsion subgroup in the first homology group of a $2$-dimensional complex on $n$ vertices can be of order $\exp(\Theta(n^2))$. 
From the data collected for $n = 50$, $60$, $75$, $100$, and $125$ we summarize the average size of $\log(|LT(n)|)$ across all 10,000 trials 
for each value of $n$ and give the results in Table \ref{tbl:size_torsion_grp}. In this table and in all future tables, values of the form $\mu \pm \sigma$ 
refer to the empirical mean $\mu$ and the standard deviation $\sigma$ of the statistic labeled by the column heading. 

\begin{table}[H]
\small\centering
\defaultaddspace=0.2em
\caption{Order of $\log(|LT(n)|)$ for $10{,}000$ trials for each value of $n$.}\label{tbl:size_torsion_grp}
\begin{tabular}{@{}r@{\hspace{12mm}}r@{}}
\\\toprule
 \addlinespace
 \addlinespace
    $n$   &  $\log(|LT(n)|)$    \\ \midrule
 \addlinespace
 \addlinespace
 \addlinespace
	50 & $12.4683 \pm 4.2591$ \\
	60 & $28.5229 \pm 5.1663$ \\
	75 & $64.8400 \pm 6.1673$ \\
	100 & $156.6246 \pm 15.5111$ \\
	125 & $291.0269 \pm 11.8643$ \\
 \addlinespace
 \addlinespace
\bottomrule
\end{tabular}
\end{table}

Using the data of the third column in Table \ref{tbl:size_torsion_grp} we found a best-fit quadratic function of $0.0328109n^2 - 2.0328n + 32.2885$. 
To test this model, we ran a few trials of $LT(n)$ for different values of $n$ and compared the average of $\log(|LT(n)|)$ found in these trials 
to the size predicted by the quadratic regression. As in other experiments we omitted trials where $LT(n)$ returned the trivial group. 
Table \ref{tbl:size_torsion_comp} shows the results.
\begin{table}[H]
\small\centering
\defaultaddspace=0.2em
\caption{Predicted and empirical size of $\log(|LT(n)|)$.}\label{tbl:size_torsion_comp}
\begin{tabular}{@{}r@{\hspace{4mm}}r@{\hspace{4mm}}r@{\hspace{4mm}}r@{}}
\\\toprule
 \addlinespace
 \addlinespace
    $n$   & Number of Trials  & Predicted Size of $\log(|LT(n)|)$ & Empirical Size of $\log(|LT(n)|)$   \\ \midrule
 \addlinespace
 \addlinespace
 \addlinespace
150 & 100 & 465.614 & 464.743 \\
175 & 25 & 681.382 & 685.413\\ 
200 & 2 & 938.164 & $948.902$ \\ 
250 & 1 & 1574.770 & $1590.538$ \\
260 & 1 & 1721.780 & $1742.036$ \\
270 & 1 & 1875.350 & $1889.057$\\
 \addlinespace
 \addlinespace
\bottomrule
\end{tabular}
\end{table}

Additionally, we examined when the largest torsion group appeared. In Table~\ref{tbl:appearance_torsion_grp} we summarize 1000 trials 
at each value of $n$ and record the average number of faces in the complex when the largest torsion group appeared. In order to quantify 
how close this is to the asymptotic homological phase transition, the third column, labeled $c$, gives a number determined by the number 
of vertices $n$ and the (average) number of faces $f$ in the following way:
$$c := \frac{nf}{\binom{n}{3}} = \frac{3f}{\binom{n - 1}{2}}.$$
This formula comes from the fact that the expected number of faces in $Y \sim Y_2(n, c/n)$, where $c$ is a constant, 
is $$\E[f(Y)] =  \frac{c}{n} \binom{n}{3} = \frac{c}{3} \binom{n - 1}{2},$$ by taking $\E[f(Y)]$ to be the empirical average given in the second column. 
Thus, to see that the torsion burst corresponds to the homological phase transition we look for the number in the third column to be close 
to the constant $c_2$ as defined exactly in \cite{ALLM}, and which to four decimal places is 2.7538. Therefore, we have strong evidence here 
pointing to the torsion burst occurring just before the phase transition.

\begin{table}[H]
\small\centering
\defaultaddspace=0.2em
\caption{Number of faces when the largest torsion group appeared.}\label{tbl:appearance_torsion_grp}
\begin{tabular}{@{}r@{\hspace{12mm}}r@{\hspace{12mm}}r@{}}
\\\toprule
 \addlinespace
 \addlinespace
     Vertices  &  Number of Faces  &  $c$     \\ \midrule
 \addlinespace
 \addlinespace
 \addlinespace
50 & 1061.413 $\pm$ 12.8200  & $2.70769 \pm 0.0327$\\
60 & 1548.091 $\pm$ 15.0146 & $2.71436 \pm 0.0263$ \\
70 & 2128.227 $\pm$ 18.3770 & $2.72152 \pm 0.0235$ \\
80 & 2798.862 $\pm$ 19.8927  & $2.72528 \pm 0.0193$ \\
90 & 3562.576 $\pm$ 22.4019 & $2.72925 \pm 0.0172$\\
100 & 4415.850 $\pm$ 24.6725 & $2.73089 \pm 0.0153$\\
110 & 5362.723 $\pm$ 27.4617 & $2.73330 \pm 0.0140$ \\
 \addlinespace
 \addlinespace
\bottomrule
\end{tabular}
\end{table}

\subsection{Other Cohen--Lenstra heuristics in the torsion burst}\label{sec:othergroups}

To this point, we have focused our attention only on the largest group in the torsion burst
in an instance $\Y_2(n)$ of the stochastic Linial--Meshulam process. 
However, as Table~\ref{tbl:75vertices} shows, it is also possible (and common) for other torsion groups to appear within the torsion burst. 
In this section we discuss experimental evidence and establish conjectures for Cohen--Lenstra heuristics in these other torsion groups. 

So far, we have understood the phrase ``torsion burst" to refer to the apparent torsion in homology around the time that the first nontrivial cycle appear. 
For our experiments above this understanding is sufficient. Here though, it is helpful to have a precise definition of the torsion burst in $\Y_d(n)$. 
\begin{definition}
Given an instance of $\Y_d(n)$ let $LT(n)$ refer to the largest torsion group which appears in $H_{d - 1}(Y_d(n, m))$, where ties 
between nonisomorphic torsion groups of maximum size are broken by the group which appears first in the stochastic process. 
Let $m_0$ be the first time $LT(n)$ appears in the codimension-$1$ homology group over $\Y_d(n)$. The \textit{torsion burst} of $\Y_d(n)$ 
is the unique maximal consecutive set of states $B$ in $\Y_d(n)$ so that $Y(n, m_0) \in B$ and every state in $B$ realizes torsion 
in codimension-$1$ homology. The \textit{duration} of the torsion burst is $|B|$. 
\end{definition}

Given an instance $\Y_2(n)$ which has a torsion burst, let $G_0 = LT(n)$ and let $m_0$ be the first place where $G_0$ 
appears in the first homology group. Starting from $m_0$, we successively inspect the torsion parts 
of $H_1(Y_2(n, m_0-1))$, $H_1(Y_2(n, m_0-2)), \dots$ The first group we find different from $G_0$ we denote by $G_{-1}$,
the second group we find different from the previously recorded groups $G_0$ and $G_{-1}$ we denote by $G_{-2}$, etc.,
until we get the trivial group as $G_{-l}$ for some $l$. 
For each $0<k \leq l$ we refer to  $G_{-k}$ as the $k$th \emph{subcritical torsion group}. 
For $k > l$, we say that the $k$th subcritical torsion group is undefined. 

Similarly, we can define the $k$th \emph{supercritical torsion group}. Starting again with $G_0 = LT(n)$ and $m_0$ 
the first place $LT(n)$ appears in the first homology group, successively inspecting the torsion parts 
of $H_1(Y_2(n, m_0+1))$, $H_1(Y_2(n, m_0+2)), \dots$ allows to define groups $G_{1}$, $G_{2}, \dots$, 
until the trivial group is reached as $G_{l'}$ for some~$l'$. 

Note that, on the formal level, our definition is not symmetric for the subcritical and the supercritical torsion groups,
as the search for new groups always begins at $m_0$ where $LT(n)$ appears for the first time.
In fact, in Table~\ref{tbl:75vertices} there is exactly one occurrence of $LT(n)$, but in other cases we saw $LT(n)$
persisting for a couple of steps, so choosing $m_0$ as the step where $LT(n)$ appears for the last time
would also be an option. Most instances of the torsion burst in $\Y_2(n)$ we saw in our experiments are \emph{unimodal}
in the sense that before we reach $LT(n)$ the size of the torsion group is monotone increasing (though not strictly increasing), 
and that after we reach $LT(n)$ the size of the torsion group is monotone decreasing. Moreover, for a unimodal torsion burst 
we see that the subcritical groups are iterative subgroups of $LT(n)$ and the supercritical groups are iterative quotients of $LT(n)$. 
This can be checked routinely by the Mayer--Vietoris sequence. 

It is possible that the torsion burst could be non-unimodal or there could even be cases where the torsion appears, disappears, 
and reappears again, which means that we inspect only the component that contains (the first occurrence of) $LT(n)$, 
but such cases are rare. This is why we include the condition of ``consecutive" in the definition of the torsion burst.

For a particular unimodal instance of $\Y_2(n)$, let $\mathcal{G}^-(\Y_2(n))$ denote the number of defined, nontrivial subcritical groups, 
and let $\mathcal{G}^+(\Y_2(n))$ denote the number of defined, nontrivial supercritical groups.
We will call the value $\mathcal{G}(\Y_2(n))=\mathcal{G}^-(\Y_2(n)) + \mathcal{G}^+(\Y_2(n)) + 1$ the number of phases in the torsion burst 
of $\Y_2(n)$, if the torsion burst exists and is unimodal (where we might miss cases where torsion groups reappear). 
Note that the ``$+1$'' in the formula for $\mathcal{G}(\Y_2(n))$ comes from adding one to count $LT(n)$.

Counting the number of phases 
is slightly different than counting the number of distinct torsion groups within the torsion burst. In counting just the number of phases 
we ignore the duration of any particular group, but we will count a group twice if it appears before and after $LT(n)$. 
Furthermore, we extend these definitions to non-unimodal torsion bursts. In that case the number of subcritical 
and the number of supercritical groups are counted with multiplicity in case some group isomorphism class appears 
as $G_{i}$ and $G_{j}$ for different $i$ and $j$, even if $i$ and $j$ are both positive or both negative. 
The number of phases is defined in the exact same way in either case.

\begin{example}
In the instance of $\Y_2(75)$ shown in Table \ref{tbl:75vertices}, the largest torsion group is  $\Z/2\Z \times \Z/79040679454167077902597570\Z$, 
the first subcritical torsion group is $\Z/2 \Z \times \Z/2 \Z$, the second subcritical torsion group is $\Z / 2\Z$, the third subcritical torsion group 
is the trivial group, and for every $k \geq 4$, the $k$th subcritical group is undefined. On the other side, the first supercritical torsion group 
is $\Z / 2\Z$, the second supercritical torsion group is trivial, and all higher supercritical torsion groups are undefined. The number of phases 
in this instance of $\Y_2(75)$ is four. We observe that the number of phases is different than the number of distinct groups since $\Z/2\Z$ 
is both the first supercritical torsion group and the second subcritical torsion group. 
\end{example}

To examine the smaller torsion groups in homology, we ran 10,000 trials at $n = 60$ and recorded the entire torsion burst. 
Based on the data collected we make the following conjectures.
\begin{conjecture}\label{cj:subcritical}
For $k \in \N$, let $\lambda_k$ denote the probability distribution on the set of all abelian groups which assigns probability proportional 
to $\dfrac{1}{|G|^k|\Aut(G)|}$ to each finite abelian group $G$. Then for each $k$, the conditional distribution on the $k$th subcritical torsion group 
of $\Y_2(n)$ given that it exists converges to $\lambda_k$.
\end{conjecture}

\begin{conjecture}\label{cj:supercritical}
Let $\lambda_k$ be as above. For each $k$, the conditional distribution on the $k$th supercritical torsion group of $\Y_2(n)$ given that it exists 
converges to $\lambda_k$.
\end{conjecture}

Note that for these conjectures it is not necessary to restrict to the $q$-part of the subcritical or supercritical group. 
The distribution which we claim is known to be a probability distribution on the set of all finite abelian groups. 
Indeed, this distribution appears for example in \cite{CL, Koplewitz, MatchettWood, MatchettWood2}. 
For any $k \geq 1$, the constant of proportionality for $\lambda_k$ is known to be 
$\prod_{p \text{ prime}}\prod_{i = k + 1}^{\infty} (1 - 1/p^{-i}) = \prod_{i = k+ 1}^{\infty} \zeta(i)^{-1} < \infty$, 
where $\zeta$ denotes the Riemann zeta function. 

Observe the symmetry between the two conjectures. For each $k$, we conjecture the same limiting conditional distribution 
for $k$-subcritical and for $k$-supercritical torsion groups. There is no immediately obvious reason for this symmetry. 
Indeed, as we pointed out above, in the unimodal case, the subcritical torsion groups are iterated subgroups of $LT(n)$ 
and the supercritical torsion groups are iterated quotients, so a priori they could behave quite differently. 
Nevertheless, our conjectures point to a remarkable symmetry.

 Tables \ref{tbl:smallgroups1}, \ref{tbl:smallgroups2}, and \ref{tbl:smallgroups3} summarize the results of 10,000 trials at $n = 60$ 
 and compares them to the conjectured Cohen--Lenstra distributions by providing the ratios of the number of instances a particular 
 sub- or supercritical group was trivial to the number of instances where it was a particular group, for common groups, 
 as in Tables \ref{tbl:2groups}, \ref{tbl:3groups}, and \ref{tbl:5groups}. For each $k$, we condition on the event 
 that the $k$th sub- or the $k$th supercritical group exists, the data in the tables reflects this. 
Now, the differences between the observed and conjectured ratios are larger in the Tables \ref{tbl:smallgroups1}, \ref{tbl:smallgroups2}, and \ref{tbl:smallgroups3},  
especially in Table \ref{tbl:smallgroups2} and \ref{tbl:smallgroups3}, than in the Tables \ref{tbl:2groups}, \ref{tbl:3groups}, and \ref{tbl:5groups}, but
this can be explained by the sample size and the high concentration of instances of the trivial group. As further evidence for our conjecture, 
we found the total variation distance between the empirical distribution and the conjectured distribution for the subcritical and supercritical groups
to be less than 0.06 in all cases.

\begin{table}
\small\centering
\defaultaddspace=0.2em
\caption{Observed ratios for the first subcritical torsion group and the first supercritical torsion group for 10,000 instances of $\Y_2(60)$, 
and a comparison to the conjectured limiting ratios.}\label{tbl:smallgroups1}
\begin{tabular*}{\linewidth}{@{}l@{\extracolsep{\fill}}r@{\extracolsep{\fill}}r@{\extracolsep{\fill}}r@{}}
\\\toprule
 \addlinespace
 \addlinespace
    Group & 1st Subcritical Group & 1st Supercritical Group & Cohen--Lenstra Ratio    \\ \midrule
 \addlinespace
 \addlinespace
 \addlinespace
Trivial Group & 1 & 1 & 1 \\
$\Z/2\Z$ & 1.9482 & 2.0793 & 2 \\
$\Z/3\Z$ & 5.7533 & 6.4367 & 6 \\
$\Z/4\Z$ & 7.7613 & 8.0608 & 8 \\
$\Z/6\Z$ & 11.9222 & 11.3984 & 12 \\
$\Z / 5 \Z$ & 18.1864 & 21.1449 & 20 \\
$\Z/2\Z \times \Z/2 \Z$ & 23.3261 & 21.9950 & 24 \\
$\Z / 8 \Z$ & 35.7667 & 30.3958 & 32 \\
$\Z/10 \Z$ & 36.3729 & 36.7815 & 40 \\
$\Z / 7 \Z$ & 42.4950 & 38.0609 & 42 \\
$\Z/12 \Z$ & 49.3333 & 42.0865 & 48 \\
$\Z / 9 \Z$ & 54.3291 & 54.0370 & 54 \\
$\Z / 2 \Z \times \Z / 4 \Z$ & 58.7945 & 61.6479 & 64 \\
$\Z / 14 \Z$ & 93.3043 & 91.1875 & 84 \\
 \addlinespace
 \addlinespace
\bottomrule
\end{tabular*}
\end{table}

\begin{table}
\small\centering
\defaultaddspace=0.2em
\caption{Observed ratios for the second subcritical torsion group and the second supercritical torsion group for 10,000 instances 
of $\Y_2(60)$, and a comparison to the conjectured limiting ratios; there were 5,708 trials where the second subcritical torsion group 
was defined and 5,623 trials where the second supercritical torsion group was defined.}\label{tbl:smallgroups2}
\begin{tabular*}{\linewidth}{@{}l@{\extracolsep{\fill}}r@{\extracolsep{\fill}}r@{\extracolsep{\fill}}r@{}}
\\\toprule
 \addlinespace
 \addlinespace
    Group & 2nd Subcritical Group & 2nd Supercritical Group & Cohen--Lenstra Ratio    \\ \midrule
 \addlinespace
 \addlinespace
 \addlinespace
Trivial Group & 1 & 1 & 1 \\
$\Z/2\Z$ & 5.2907 & 4.7272 & 4 \\
$\Z/3\Z$ & 26.3274 & 23.5833 & 18 \\
$\Z/4\Z$ & 39.4911 & 39.3056 & 32 \\
$\Z/6 \Z$ & 107.8780 & 108.8462 & 72 \\
$\Z/2 \Z \times \Z/2 \Z$ & 126.3714 & 92.2826 & 96 \\
$\Z/5\Z$ & 152.5172 & 132.6563 & 100 \\
 \addlinespace
 \addlinespace
\bottomrule
\end{tabular*}
\end{table}

\begin{table}
\small\centering
\defaultaddspace=0.2em
\caption{Observed ratios for the third subcritical torsion group and the third supercritical torsion group for 10,000 instances of $\Y_2(60)$ 
and a comparison to the conjectured limiting ratios; there were 1,285 trials where the third subcritical torsion group was defined 
and 1,378 trials where the third supercritical torsion group was defined.}\label{tbl:smallgroups3}
\begin{tabular*}{\linewidth}{@{}l@{\extracolsep{\fill}}r@{\extracolsep{\fill}}r@{\extracolsep{\fill}}r@{}}
\\\toprule
 \addlinespace
 \addlinespace
    Group & 3rd Subcritical Group & 3rd Supercritical Group & Cohen--Lenstra Ratio    \\ \midrule
 \addlinespace
 \addlinespace
 \addlinespace
Trivial Group & 1 & 1 & 1 \\
$\Z/2\Z$ & 11.0866 & 10.9206 & 8 \\
$\Z/3 \Z$ & 104.8182 & 77.1250 & 54 \\
 \addlinespace
 \addlinespace
\bottomrule
\end{tabular*}
\end{table}

Finally, a related statistic is the duration of the torsion burst.  Before presenting the experimental data, we note that Conjectures \ref{cj:subcritical} 
and \ref{cj:supercritical} imply something about the duration of the torsion burst. Namely Conjecture \ref{cj:subcritical} (resp. Conjecture \ref{cj:supercritical}) 
implies that there is a positive probability the $k$th subcritical (supercritical) torsion group is defined for any fixed $k$, 
provided there is a positive probability there is a torsion burst. This is straightforward to compute. Let $p_0$ denote the asymptotic probability 
that $\Y_2(n)$ has a torsion burst (assuming such a probability exists). For $k \geq 1$, let $p_k$ denote the asymptotic probability 
that the $k$th subcritical torsion group of $\Y_2(n)$ is nontrivial conditioned on the event that it exists. By Conjecture~\ref{cj:subcritical}, 
$p_k > 0$ for all $k \geq 1$, in fact, $p_k = 1 - \prod_{i = k + 1}^{\infty} \zeta(i)^{-1}$. Let $q_k$ denote the asymptotic probability 
that $\Y_2(n)$ has a $k$th subcritical torsion group. Clearly, $\Y_2(n)$ has a $k$th subcritical torsion group if and only if it has a nontrivial 
$(k-1)$st subcritical torsion group. Thus we have the following recurrence $q_k = p_{k - 1}\cdot q_{k - 1}$, and $q_1 = p_0$ since a torsion burst 
implies that there is a nontrivial largest torsion group, so there is a defined first subcritical torsion group.  Thus, if $p_0 > 0$ then $q_k > 0$ 
for every positive integer~$k$. 

Moreover, by linearity of expectation the asymptotic expected number of positive integers $k$ so that $\Y_2(n)$ has a $k$th subcritical torsion group 
is given by 
$$\sum_{i = 1}^{\infty} q_i = \sum_{i = 1}^{\infty} \prod_{j = 1}^{i - 1} p_j.$$
Assuming Conjecture \ref{cj:subcritical}, this sum converges since $p_k \rightarrow 0$ as $k \rightarrow \infty$.  The same would hold 
for the supercritical torsion groups as well. Thus, we may compute the expected number of phases in the torsion burst. 
Since $\mathcal{G}^-(\Y_2(n))$ does not count the trivial subcritical torsion group we have 
$$\lim_{n \rightarrow \infty} \E[\mathcal{G}^-(\Y_2(n))] = \left(\sum_{i = 1}^{\infty} q_i\right) - 1.$$

The same expectation would hold for $\mathcal{G}^+(\Y_2(n))$. Thus the expected number of phases, assuming 
Conjectures \ref{cj:subcritical} and \ref{cj:supercritical} and that the torsion burst occurs with high probability, is asymptotically given by:

$$\lim_{n \rightarrow \infty} \E[\mathcal{G}(\Y_2(n))] = 2\left(\sum_{i = 1}^{\infty} q_i\right) - 2 + 1 = 2\left(\sum_{i = 1}^{\infty} \prod_{j = 1}^{i - 1} \left(1 - \prod_{k = j + 1}^{\infty} \zeta(k)^{-1}\right)\right) - 1.$$
A reasonable approximation for this value is 2.49524.

To collect data on the duration of the torsion burst experimentally we ran 1000 trials at $n = 50$, 60, 70, 80, 90, and 100. 
Table \ref{tbl:burstlengths} summarizes the results. 
Additionally, we should compare the number of phases in the torsion burst with the number predicted by the Cohen--Lenstra heuristics given above. 

We notice that the average duration of the torsion burst decreases with $n$. This, together with the conjectured asymptotic number of phases, 
suggests that the duration in $\Y_2(n)$ for which a particular nontrivial torsion group in the first homology group persists depends on $n$. 
There are a few trivial conditions which imply that the torsion part of $H_1(Y_2(n,m))$ is isomorphic to the torsion part of $H_1(Y_2(n, m + 1))$, 
which depend on $n$. For example, if the $(m+1)$st face completes the boundary of a tetrahedron then it will not change the torsion in homology. 
Similarly, if an edge of the $(m + 1)$st face was isolated in $Y(n,m)$ then the torsion in homology would not change. 
In the density regime we are interested in, the number of possible faces which complete a tetrahedron is expected to be $O(n)$, 
and the expected number of isolated edges is $O(n^2)$. Thus, at any stage in this regime the probability that we add a face 
which completes a tetrahedron is $O(n^{-2})$, and the probability we add a face which covers an isolated edge is $O(n^{-1})$, 
so both of these probabilities go to zero with $n$. This gives a partial explanation for why a particular nontrivial torsion group 
is asymptotically unlikely to persist for more than one state, but the full picture remains unclear.

\begin{table}
\small\centering
\defaultaddspace=0.2em
\caption{Durations of the torsion burst (1000 trials for each $n$).}\label{tbl:burstlengths}
\begin{tabular}{@{}r@{\hspace{12mm}}r@{\hspace{12mm}}r@{}}
\\\toprule
 \addlinespace
 \addlinespace
    Vertices  &   Duration & Number of Phases   \\ \midrule
 \addlinespace
 \addlinespace
 \addlinespace
50 & $6.882 \pm 6.087$ & $2.509 \pm 1.215$ \\
60 & $6.019 \pm 6.494$ & $2.438 \pm 1.212$\\
70 & $5.788 \pm 6.856$ & $2.455 \pm 1.213$\\
80 & $5.229 \pm 4.533$ & $2.398 \pm 1.205$\\
90 & $5.094 \pm 4.000$ & $2.394 \pm 1.146$\\
100 & $5.205 \pm 6.710$ & $2.440 \pm 1.173$ \\
110 & $5.097 \pm 4.232$ & $2.426 \pm 1.176$\\
 \addlinespace
 \addlinespace
\bottomrule
\end{tabular}
\end{table}

\section{Torsion in higher-dimensional complexes}
\label{sec:higher_dim}

So far all of our experiments have been in the Linial--Meshulam model with $d = 2$. The next step would be to see what happens 
in higher dimensions. It is not immediately obvious that the same torsion burst should occur in higher dimension. For one, 
the $2$-dimensional case has a nontrivial fundamental group, and so the torsion in the first homology group could be coming from 
the fundamental group. In higher dimensions, the random complex is simply-connected, and so one might expect that 
something completely different happens here.

On the other hand, as \cite{LP} points out, $Y_d(n, c/n)$ has a phase transition 
with regard to Euler characteristic at $c = c_d$ in the sense that when $c < c_d$ one has that the pure part of the complex $X$ 
obtained by collapsing $Y_d(n, c/n)$ as far as possible has negative Euler characteristic, whereas for $c > c_d$ this pure part of the complex obtained by 
performing all possible collapses  
has positive Euler characteristic. So perhaps torsion appears at the moment the Euler characteristic of this ``essential core" reaches one, 
because at that point we have something that is close to a $\Q$-acyclic complex, which, as we mentioned above, 
typically have large torsion in homology, by a result of \cite{Kalai}. In this case we might expect that the torsion burst 
is a general phenomenon in higher dimensions, and early evidence suggests this is correct. 
Tables \ref{tbl:25vertices3d}, \ref{tbl:17vertices4d}, and \ref{tbl:16vertices5d} show results of a single trial each in 3, 4, and 5 dimensions, respectively.

\begin{table}[H]
\small\centering
\defaultaddspace=0.2em
\caption{Homology groups of one instance of $\Y_3(25)$.}\label{tbl:25vertices3d}
\begin{tabular}{@{}r@{\hspace{12mm}}l@{\hspace{12mm}}l@{}}
\\\toprule
 \addlinespace
 \addlinespace
    Faces  & $H_3$  & $H_2$    \\ \midrule
 \addlinespace
 \addlinespace
 \addlinespace
1949 & $\Z^4$ & $\Z^{79}$ \\
1950 & $\Z^4$ & $\Z^{78} \times \Z/6\Z$ \\
1951 & $\Z^4$ & $\Z^{77} \times \Z/7780167918307023583785903521760\Z$ \\
1952 & $\Z^5$ & $\Z^{77} \times \Z/5\Z$  \\
1953 & $\Z^6$ & $\Z^{77}$ \\
 \addlinespace
 \addlinespace
\bottomrule
\end{tabular}
\end{table}

\begin{table}[H]
\small\centering
\defaultaddspace=0.2em
\caption{Homology groups of one instance of $\Y_4(17)$.}\label{tbl:17vertices4d}
\begin{tabular}{@{}r@{\hspace{12mm}}l@{\hspace{12mm}}l@{}}
\\\toprule
 \addlinespace
 \addlinespace
    Faces  & $H_4$  & $H_3$    \\ \midrule
 \addlinespace
 \addlinespace
 \addlinespace
1787 & $\Z^{10}$ & $\Z^{43}$ \\
1788 & $\Z^{10}$ & $\Z^{42} \times \Z/2\Z$ \\
1789 & $\Z^{10}$ & $\Z^{41} \times \Z/2\Z$ \\
1790 & $\Z^{10}$ & $\Z^{40} \times \Z/2\Z$ \\
1791 & $\Z^{10}$ & $\Z^{39} \times \Z/49234986784469188898774\Z$ \\
1792 & $\Z^{11}$ & $\Z^{39}$ \\
 \addlinespace
 \addlinespace
\bottomrule
\end{tabular}
\end{table}

\begin{table}[H]
\small\centering
\defaultaddspace=0.2em
\caption{Homology groups of one instance of $\Y_5(16)$.}\label{tbl:16vertices5d}
\begin{tabular}{@{}r@{\hspace{12mm}}l@{\hspace{12mm}}l@{}}
\\\toprule
 \addlinespace
 \addlinespace
    Faces  & $H_5$  & $H_4$    \\ \midrule
 \addlinespace
 \addlinespace
 \addlinespace
2972 & $\Z^6$ & $\Z^{37}$ \\
2973 & $\Z^6$ & $\Z^{36} \times \Z/1147712621067945810235354141226409657574376675\Z$ \\
2974 & $\Z^7$ & $\Z^{36}$ \\
 \addlinespace
 \addlinespace
\bottomrule
\end{tabular}
\end{table}

We point out here, using the results of \cite{AL} and  \cite{LP}, that for the Tables \ref{tbl:25vertices3d}, \ref{tbl:17vertices4d}, and~\ref{tbl:16vertices5d}
the expected number of faces when homology emerges  are 1,976.94, 1,805.44, and 2,992.99, respectively. 
So it still seems probable in higher dimensions that torsion in codimension-$1$ homology occurs immediately before the emergence 
of nontrivial top-dimensional homology.

In addition to the Linial--Meshulam model we also tried to compute torsion in homology of the stochastic process on random clique complexes. 
Recall that the clique complex of a graph $G$ is the largest simplicial complex having $G$ as its $1$-skeleton, i.e., 
the $d$-dimensional simplices are exactly the $d$-cliques in $G$. In the stochastic model one edge is added at a time and homology 
of the resulting clique complex is computed at each step. However, in running this process on up to 100 vertices several times 
we never found any torsion in the clique complex model. Perhaps this is not too surprising. In the Linial--Meshulam model,
torsion seems to vanish very quickly.

This suggests that many of the empty triangles are cycles which have finite order in the first homology group, 
so a randomly selected one is likely to drop the torsion. In contrast, in the clique complex model, there are no empty triangles to be torsion cycles. 
Additionally, there are $\binom{n}{d + 1}$ steps in the Linial--Meshulam stochastic model and the torsion burst seems to only last for a few faces, 
but the clique complex stochastic model only has $\binom{n}{2}$ steps, so it could be in some sense too coarse to detect torsion.


\section{Random $2$-trees}
\label{sec:Q-acyclic}
In this section, we  study uniform random $2$-trees in the $2$-skeleton of the simplex on $n$ vertices. These are $2$-dimensional analogues of uniform spanning trees on complete graphs on $n$ vertices.

A $2$-dimensional simplicial complex $T$ is a \emph{$2$-tree} if
 $T$ has a complete $1$-skeleton,
 $H_1(T, \Q) = 0$, and
 $H_2(T, \Q) = 0$.

$2$-trees are also called $\Q$-acyclic complexes. They were first studied topologically and combinatorially by Kalai \cite{Kalai}, who showed as a corollary of his higher-dimensional Cayley's theorem that the expected order of the torsion part of $H_1( T, \Z)$ for a uniform random $2$-tree on $n$ vertices is enormous. It is at least $\exp(cn^2)$ for some constant $c >0$. 

Little else seems to be known about the topology of random $2$-trees. It is apparently an open problem even to show that, with high probability, these complexes are not contractible. Our experiments suggest that there is almost always nontrivial torsion in homology, and that the torsion group is Cohen--Lenstra distributed.


\subsection{Markov chains}

We review a few basic facts and definitions about Markov chains.  Consider a discrete Markov chain with state space $A$ 
and transition matrix $P$.  We call the Markov chain \emph{irreducible} if for all pairs of states $a_1,a_2\in A$, 
there exists an integer $k\geq0$ such that $P^k_{a_1,a_2}>0$. That is, the Markov chain is irreducible if it is possible 
to transition from any state to any other state in finite time.  A state $a\in A$ is \emph{lazy} provided $P_{a,a}>0$, 
that is, if there is positive probability of not transitioning from $a$ when the current state is $a$.  
The Markov chain is called \emph{regular} provided there exists $k'>0$ such that for all $a_1,a_2\in A$, $P^{k'}_{a_1,a_2}>0$ 
(note the difference with the definition of irreducible).  It is easy to see that a finite, irreducible Markov chain with at least one lazy state is regular.  
A distribution $\pi$ on the state space $A$ satisfies \emph{detailed balance} provided for all pairs of states 
$a_1,a_2\in A$, $P_{a_1,a_2}\pi_{a_1} = P_{a_2,a_1}\pi_{a_2}$.  It is well known (see, for example \cite{Levin:2009}) 
that given a finite, regular Markov chain with state space $A$, transition matrix $P$, and distribution $\pi$ on $A$ satisfying detailed balance, 
$\pi$ is the unique limit distribution of the Markov chain given any initial distribution.  That is, given any distribution $\tilde{\pi}$ on $A$,
\begin{displaymath}
\lim_{n\to\infty}P^n\tilde{\pi} = \pi.
\end{displaymath}
We will use these basic facts about Markov chains to demonstrate that the chain we build 
on the set $\mathcal{T}^2(n)$ of all $2$-dimensional $\Q$-acyclic complexes on $n$ vertices 
approximately samples from the uniform distribution on $\mathcal{T}^2(n)$.

We implement a Markov chain on the state space $\mathcal{T}^2(n)$ in the following way.  Given $T\in\mathcal{T}^2(n)$, 
uniformly choose $\sigma\in f_2(T)$ and $\tau\in f_2(K_n^2)\setminus f_2(T)$.  Let $T'$ be the subcomplex of $K_n^2$ 
with all the same cells as $T$, except $f_2(T') = (f_2(T)\setminus\{\sigma\})\cup\{\tau\}$.  If $T'\in\mathcal{T}^2(n)$ then the Markov chain 
transitions to $T'$ otherwise it remains at $T$.  Call the transition matrix for this Markov chain~$P$.  
Note that $\lvert f_2(T)\rvert = \binom{n-1}{2}$ and $\lvert f_2(K_n^2)\setminus f_2(T)\rvert = \binom{n}{3}-\binom{n-1}{2} = \binom{n-1}{3}$.  
We see then that $P$ is given by

\begin{displaymath}
P_{T_1,T_2} = 
\begin{cases}
\frac{1}{\binom{n-1}{2}\binom{n-1}{3}}, & \lvert f_2(T_1)\bigtriangleup f_2(T_2)\rvert = 2, \\
0, & \lvert f_2(T_1)\bigtriangleup f_2(T_2)\rvert > 2, \\
1-\sum_{T\neq T_1}{P_{T_1,T}}, & T_1 = T_2. \\
\end{cases}
\end{displaymath}
Since $P_{T_1,T_2} = P_{T_2,T_1}$ we know that the uniform distribution on $\mathcal{T}^2(n)$ satisfies detailed balance 
and further since $\mathcal{T}^2(n)$ is the collection of maximal independent sets of a matroid the exchange principle shows 
the Markov chain is irreducible.  Note finally that there are subcomplexes $Y$ of $K_n^2$ with $\lvert f_2(Y)\rvert = \binom{n-1}{2}$ 
that are not in $\mathcal{T}^2(n)$.  For example, any subcomplex of $K_n^2$ with complete $1$-skeleton that contains within its $2$-faces 
a $2$-cycle such as the boundary of the tetrahedron $\{1,2,3\},\{1,2,4\},\{1,3,4\}$, and $\{2,3,4\}$.  This gives us that there is at least one 
(in reality many more than one) $T\in\mathcal{T}^2(n)$ with $P_{T,T}>0$.  Thus the Markov chain has states that are lazy with positive probability 
and so is regular.  Given that it is irreducible, regular, and the uniform distribution satisfies detailed balance, we can conclude 
that the uniform distribution is the unique limiting distribution for any initial distribution.  Thus our Markov chain will sample approximately 
from the uniform distribution on $\mathcal{T}^2(n)$.

Sampling from the uniform distribution on $\mathcal{T}^2(n)$ and then applying $H_1(-,\mathbb{Z})$ yields a distribution $\pi(n)$ on finite abelian groups. Taking finitely many samples from the Markov chain on $\mathcal{T}^2(n)$ and applying $H_1(-,\mathbb{Z})$ gives an empirical distribution 
$\tilde{\pi}(n)$ that approximates $\pi(n)$.  We will analyze if $\pi(n)$ also appears to conform to the Cohen--Lenstra heuristics. 
Mirroring our work on the Linial--Meshulam model we look only at the $p$-part of our distribution. That is, for every prime $p$ 
we have a distribution on $p$-groups, $\pi^p(n)$, given by sampling a group from $\pi(n)$ and looking at its Sylow $p$-subgroup. 
The corresponding empirical distributions will be denoted $\tilde{\pi}^p(n)$. 

\subsection{Experiments and Results}
%

{Apparently, the question of whether the bases-exchange Markov chain mixes quickly for every matroid is a well-known open question in probability \cite{Peres18}. The Markov chain mixes quickly for ``balanced'' matroids \cite{FM92} where matroid elements are negatively correlated, but unfortunately the $2$-tree matroid we study is not balanced. On 6 vertices, for example, there are 46620 2-dimensional $\Q$-acyclic complexes. Half of these contain the face $[1, 2, 3]$, and half contain the face $[4, 5, 6]$. However, 11664 contain both. So the probability that a uniform random $2$-tree on 6 vertices contains $[1, 2, 3]$ and $[4, 5, 6]$ is $$11664/46620 > 1/4,$$
contradicting negative correlation.

We ran the Markov chain until time $2t_0$ where $t_0$ was the first time that the edge degrees was a set of consecutive integers without gaps. This was a somewhat arbitrary choice, but based on the intuition that the degrees would all get closer to concentrated around their expected value. Even though we are not aware of any rigorous results on the mixing time of our Markov chain, our experiments makes us guess that at least the torsion may already be close to its limiting distribution by this point.

We conjecture the same Cohen--Lenstra distribution for torsion in the random $2$-tree that we conjectured for the torsion burst in the Linial--Meshulam model.

\begin{conjecture}
For a fixed prime $q$, the Sylow $q$-subgroup of $H_1(X)$, where $X$ is drawn uniformly from $\mathcal{T}^2(n)$, 
is asymptotically distributed according to the Cohen--Lenstra heuristic which assigns probability 
$\dfrac{\prod_{k = 1}^{\infty} (1 - q^{-k})}{|\Aut(G)|}$ to any $q$-group $G$.
\end{conjecture}

The resulting homology groups are torsion groups by construction, and even though we conjecture the same limiting distribution, the groups we found for random $2$-trees were somewhat larger than the torsion groups we found in the Linial--Meshulam process. For example, the average order of the torsion group for our experiments on 100 vertices is $9.92\cdot 10^{96}$, with a standard deviation of $3.99\cdot 10^{98}$.

We generated $\Q$-acyclic complexes on 50, 60, 75, and 100 vertices as described above.  
Then for small primes $q$ we calculated the $q$-part of the first homology of the complexes.  For a few of the more common groups 
Tables \ref{tbl:2trees}, \ref{tbl:3trees}, and \ref{tbl:5trees} compare for each $p$-group $G$ the ratio of instances 
where the Sylow $q$-subgroup was isomorphic to $G$ to the number of instances where the Sylow $q$-subgroup was trivial, 
as in Tables \ref{tbl:2groups}, \ref{tbl:3groups}, and \ref{tbl:5groups}. For each value of $n$ we generated 10,000 $\Q$-acyclic complexes.

\begin{table}
\small\centering
\defaultaddspace=0.2em
\caption{$2$-subgroups of the torsion in random $\Q$-acyclic complexes}\label{tbl:2trees}
\begin{tabular*}{\linewidth}{@{}l@{\extracolsep{\fill}}r@{\extracolsep{\fill}}r@{\extracolsep{\fill}}r@{\extracolsep{\fill}}r@{\extracolsep{\fill}}r@{}}
\\\toprule
 \addlinespace
 \addlinespace
    Group & $n = 50$ & $n = 60$ & $n = 75$ & $n = 100$ & Conjectured Ratio \\\midrule
 \addlinespace
 \addlinespace
 \addlinespace
No $2$-torsion & 1& 1& 1& 1& 1 \\
$\mathbb{Z}/2\mathbb{Z}$& 0.938748 & 1.024867 & 0.986878 & 1.020466 & 1 \\
$\mathbb{Z}/4\mathbb{Z}$& 1.909153 & 1.997922 & 2.000000 & 1.893910 & 2 \\
$\mathbb{Z}/8\mathbb{Z}$& 3.879720 & 3.841545 & 4.082857 & 3.881879 & 4 \\
$\mathbb{Z}/2\mathbb{Z}\times\mathbb{Z}/2\mathbb{Z}$& 5.503968 & 5.747012 & 5.614931 & 6.426667 & 6 \\
$\mathbb{Z}/16\mathbb{Z}$& 7.397333 & 7.904110 & 8.073446 & 8.406977 & 8 \\
$\mathbb{Z}/4\mathbb{Z}\times\mathbb{Z}/2\mathbb{Z}$& 7.149485 & 8.266476 & 7.916898 & 7.630607 & 8 \\
$\mathbb{Z}/32\mathbb{Z}$& 14.373057 & 15.104712 & 14.148515 & 14.830769 & 16 \\
$\mathbb{Z}/8\mathbb{Z}\times\mathbb{Z}/2\mathbb{Z}$& 15.158470 & 16.485714 & 16.713450 & 15.221053 & 16 \\
 \addlinespace
 \addlinespace
\bottomrule
\end{tabular*}
\end{table}

\begin{table}
\small\centering
\defaultaddspace=0.2em
\caption{$3$-subgroups of the torsion in random $\Q$-acyclic complexes}\label{tbl:3trees}
\begin{tabular*}{\linewidth}{@{}l@{\extracolsep{\fill}}r@{\extracolsep{\fill}}r@{\extracolsep{\fill}}r@{\extracolsep{\fill}}r@{\extracolsep{\fill}}r@{}}
\\\toprule
 \addlinespace
 \addlinespace
    Group & $n = 50$ & $n = 60$ & $n = 75$ & $n = 100$ & Conjectured Ratio \\\midrule
 \addlinespace
 \addlinespace
 \addlinespace
No $3$-torsion & 1& 1& 1& 1& 1 \\
$\mathbb{Z}/3\mathbb{Z}$& 2.030946 & 1.965724 & 1.987892 & 1.940909 & 2 \\
$\mathbb{Z}/9\mathbb{Z}$& 6.222712 & 5.988159 & 5.900634 & 5.924226 & 6 \\
$\mathbb{Z}/27\mathbb{Z}$& 18.565789 & 17.116923 & 17.443750 & 17.678344 & 18 \\
$\mathbb{Z}/3\mathbb{Z}\times\mathbb{Z}/3\mathbb{Z}$& 52.259259 & 47.956897 & 40.158273 & 46.647059 & 48 \\
$\mathbb{Z}/81\mathbb{Z}$& 49.946903 & 58.557895 & 64.160920 & 52.866667 & 54 \\
 \addlinespace
 \addlinespace
\bottomrule
\end{tabular*}
\end{table}

\begin{table}
\small\centering
\defaultaddspace=0.2em
\caption{$5$-subgroups of the torsion in random $\Q$-acyclic complexes}\label{tbl:5trees}
\begin{tabular*}{\linewidth}{@{}l@{\extracolsep{\fill}}r@{\extracolsep{\fill}}r@{\extracolsep{\fill}}r@{\extracolsep{\fill}}r@{\extracolsep{\fill}}r@{}}
\\\toprule
 \addlinespace
 \addlinespace
    Group & $n = 50$ & $n = 60$ & $n = 75$ & $n = 100$ & Conjectured Ratio \\\midrule
 \addlinespace
 \addlinespace
 \addlinespace
No $5$-torsion & 1& 1& 1& 1& 1 \\
$\mathbb{Z}/5\mathbb{Z}$& 4.068060 & 3.874423 & 4.242123 & 4.186783 & 4 \\
$\mathbb{Z}/25\mathbb{Z}$& 18.514634 & 21.055710 & 18.358852 & 20.334218 & 20 \\
 \addlinespace
 \addlinespace
\bottomrule
\end{tabular*}
\end{table}

We estimated the total variation distances between the empirical distributions coming from our data sets and the conjectured distributions for primes $q=2, 3, 5, \dots, 23$; for $n=100$, all these distances were less than $0.04$.

\section{A hitting time conjecture in the Linial--Meshulam model} \label{sec:hit}

At this point we briefly review what is known about the Linial--Meshulam model in the $p=c/n$ regime. 
In \cite{AL} and \cite{LP}, a one-sided sharp threshold of $p=c_d/n$ is established for the property that $H_d( Y_d(n, p)) \neq 0$, 
with $c_d$ a constant depending only on $d$ and explicitly given in \cite{AL}. We do not give the precise definition of $c_d$ here, 
but $c_2 \approx 2.75383$.  
  Beyond this, however, \cite{LP} says more about how the complex changes at $c_d/n$. We have the following theorems:
\begin{theorem}\cite[Theorem 1.1]{LP}
For $c < c_d$, with high probability, $H_d(Y)$ is generated by a bounded number of copies of the boundary of the $(d + 1)$-simplex. 
\end{theorem}
\begin{theorem}\cite[Theorem 1.3]{LP}
For $c > c_d$, let $t_c$ be the smallest positive root of $t = e^{-c(1 - t)^d}$. Denote $\beta_d =\dim H_d(Y_d(n, c/n), \R)$.  With high probability,
$$\beta_d = \left(1 + o(1) \right) \left( \binom{n}{d} \left(ct_c(1 - t_c)^d + \dfrac{c}{d + 1}(1 - t_c)^{d + 1} - (1 - t_c)\right) \right).$$
\end{theorem}

Furthermore, not only do we know that a phase transition occurs in $H_d$ at $c_d/n$, the following result about the shadow is also proved by \cite{LP}. 
The \emph{shadow} of a $d$-dimensional simplicial complex $X$ with complete $(d-1)$-skeleton over a field $\F$ is denoted $SH_{\F}(Y)$ 
and is defined first by \cite{LNPR} as
$$SH_{\F}(Y) = \{ f \in \Delta_{n -1}^{(d)}\backslash f_d(Y) : \beta_d(Y \cup \{f\}; \F) > \beta_d(Y, \F)\}.$$
\begin{theorem}\cite[Theorem 1.4]{LP}
Let $Y \sim Y_d(n, c/n)$ for $d \geq 2$, $c > 0$ fixed. If $c < c_d$ then with high probability $|SH_{\R}(Y)| = \Theta(n)$, 
and if $c > c_d$ then $|SH_{\R}(Y)| = \Theta(n^{d + 1})$.
\end{theorem}

In \cite{LP} this result about the shadow is interpreted as a higher-dimensional analogue of the giant component in the Erd\H{o}s--R\'{e}nyi random graph. 
Recall that in the Erd\H{o}s--R\'{e}nyi random graph model $G(n, p)$ we have a phase transition at $p = 1/n$. 
Namely for $c < 1$ and $G \sim G(n, c/n)$ one has that with high probability every component of $G$ is on $O(\log n)$ vertices 
and that for $c > 1$ and $G \sim G(n, c/n)$ one has that with high probability $G$ has a unique ``giant component" on $\Theta(n)$ vertices 
and all other components are on $O(\log n)$ vertices; see for example \cite{JLR}. With this in mind, we make the following conjecture, 
which, if true, would enrich the analogy between the giant component of $G(n, p)$ and the giant shadow of $Y_d(n, p)$ in a way that contrasts 
some key differences between the two. As in \cite{ALLM}, a \emph{core} in a $d$-complex is a subcomplex  in which every $(d-1)$-face is contained in at least two $d$-faces. A core is an obstruction to collapsibility.

\begin{conjecture} \label{sec:hitting}
Let $\Y_d(n)$ be the stochastic Linial--Meshulam process in $d$ dimensions for $d \geq 2$. Then there exists a constant $\delta$ depending on $d$ 
so that with high probability there exists $m_0 \in \{1,\dots, \binom{n}{d + 1}\}$ so that the following three events occur.

\vspace{1mm}

\begin{enumerate}
\item The torsion part of $H_{d - 1}(Y_d(n, m_0))$ is isomorphic to $LT(n)$. \\[-2mm]
\item $Y_d(n, m_0)$ contains a core $Y'$ which spans the entire vertex set and has $H_{d - 1}(Y') \cong LT(n)$, but $Y_d(n, m_0 - 1)$ does not contain such a core. \\[-2mm]
\item $|SH_{\R}(Y_d(n, m_0))| \geq \delta n^{d+1}$, but $|SH_{\R}(Y_d(n, m_0 - 1))| \leq O(n)$.
\end{enumerate}
\end{conjecture}

%
%

Our conjecture allows us to define a higher-dimensional analogue of the giant component in terms of the faces of the complex rather 
than in terms of the shadow. Given a $d$-dimensional simplicial complex $Y$ on $n$ vertices we say that $Y$ contains a \textit{homological giant} $Y'$ 
if $Y'$ is a core of $Y$, $H_{d - 1}(Y')$ is finite, and $V(Y) = V(Y')$. The conjecture says the hitting time for the emergence of a homological giant 
is the same as the hitting time for the emergence of the giant shadow. If this is true, then with high probability the homological giant 
has its top dimensional homology generated by a bounded number of $(d+1)$-simplex boundaries at the moment it first appears. 
This means that after removal of one face from each of these $(d+1)$-simplex boundaries we have that the remaining part 
of the homological giant is a spanning simplicial tree in the sense of Duval, Klivans, and Martin \cite{Duval:2016}. That is, this is a subcomplex with $\beta_{d - 1} = \beta_{d} = 0$. This is a more general definition than the $2$-trees discussed earlier and their $d$-dimensional analogues, in the sense that the $(d-1)$-skeleton is not required to be complete.
 
Linial and Peled conjecture \cite{LP3} that in a single step in the process $\Y_d(n)$, the size of the shadow jumps 
from $O(n)$ to $\Theta(n^{d + 1})$, and they suggest the problem of trying to study the structure of this complex at the moment 
the giant shadow appears. We conjecture the same behavior as part (3) of the hitting-time conjecture (Conjecture~\ref{sec:hitting}), 
so if true, our conjecture would provide a partial answer to the problem posed by Linial and Peled in~\cite{LP3}. 
Quick experiments show evidence for the conjecture. Indeed, in 100 trials on each of $n= 50, 60, 70, 100$ we found strong evidence 
that the emergence of the giant shadow coincided with the torsion burst---these events two coincided in 99\% of the trials. 

\bigskip

\noindent {\bf Acknowledgements.} M.K.\ is grateful to Sam Payne and Yuval Peled for helpful and inspiring conversations. All four authors are grateful to ICERM for visits during the special semester program \textit{Topology in Motion} in the fall of 2016, where much of this work was completed. F.H.L.\ is also thankful to HIM Bonn for a stay for the Special Hausdorff Program \textit{Applied and Computational Algebraic Topology} in 2017.

\bibliography{ResearchBibliography}
\bibliographystyle{amsplain}

\end{document}